\documentstyle[12pt]{article}
\title{Uniform Kadec-Klee Property in Banach Lattices} 
\author{Fouad Chaatit  and  Mohamed A. Khamsi.} 

\date{}
\begin{document}
\maketitle 
\begin{abstract}
{We prove that a Banach lattice $X$ which does not contain the
$l^n_{\infty}$-uniformly has an equivalent norm which is uniformly
Kadec-Klee for a natural topology $\tau$ on $X$.   In case
the Banach lattice is purely atomic, the topology $\tau$ is
the coordinatewise convergence topology.}
\end{abstract}
\vspace*{3in}

\newpage      
\baselineskip=24pt
\font\tfont= cmbx10 scaled \magstep3
\font\afont= cmcsc10 scaled \magstep2
\centerline{\tfont Uniform Kadec-Klee Property in Banach lattices} 
\bigskip
\centerline{\afont F. Chaatit and  M. A. Khamsi }
\bigskip
\begin{abstract}
{We prove that a Banach lattice $X$ which does not contain the
$l^n_{\infty}$-uniformly has an equivalent norm which is uniformly
Kadec-Klee for a natural topology $\tau$ on $X$.   In case
the Banach lattice is purely atomic, the topology $\tau$ is
the coordinatewise convergence topology.}
\end{abstract}

\newtheorem{thm}{Theorem }[chapter]
\newtheorem{lem}[thm]{Lemma} 
\newtheorem{cor}[thm]{Corollary}
\newtheorem{prop}[thm]{Proposition}          
\newtheorem{ex}{Exercise}                             
\newenvironment{proof}{\medskip \par \noindent {\bf Proof:}\ }{\hfill $\Box$ 
                       \medskip \par}

\setcounter{section}{0}
\setcounter{thm}{0}
\setcounter{ex}{0} 
\setcounter{page}{1}

\newcommand{\bdfn}{\noindent{\bf Definition:\,} \begin{em}}
\newcommand{\edfn}{\end{em}}
\newcommand{\ntn}{\noindent{\bf Notations:\,\,}}

\newcommand{\rmk}{\noindent{\bf Remark\,\,}}


\newcommand{\del}{\mbox{$\delta$}}
\newcommand{\lftyn}{\mbox{${l_{\infty}^{n}}$}}
\newcommand{\norm}[1]{\mbox{$\|#1\|$}}             
\newcommand{\som}{\mbox{$\displaystyle\sum_{n=0}^{\infty}$}}
\newcommand{\son}{\mbox{$\displaystyle\sum_{n=1}^{\infty}$}}
\newcommand{\reals}{\mbox{\bf R}}
\newcommand{\complexes}{\mbox{\bf C}}
\newcommand{\integers}{\mbox{\bf N}}
\newcommand{\smint}{\mbox{\scriptsize {\bf N}}}
\newcommand{\complex}{\mbox{\bf C}}
\newcommand{\ep}{\mbox{$\varepsilon$}}
\newcommand{\B}{\mbox{${\cal B}_1(K)$}}
\newcommand{\pf}{\noindent{\bf Proof:\,\ }}
\newcommand{\ra}{\mbox{$\longrightarrow$}}
\def\qed{~\hfill~\blackbox\medskip}
\def\blackbox{\hbox{\vrule width6pt height7pt depth1pt}}

\newcommand{\osm}{\mbox{$(\Omega,\Sigma,\mu)$}}
\newcommand{\om}{\mbox{$\omega$}}

\section{Introduction} 

The uniform Kadec-Klee (in
short UKK) property was introduced by R. Huff in \cite{bib:H}.  Banach spaces with
Schur property, as well as uniformly convex spaces, enjoy the UKK
property. P. Enflo characterized Banach spaces which have an equivalent
 uniformly convex  norm as those which are super-reflexive.  Similarly
R. Huff raised the following question :\\
{\em Under what conditions does a Banach space possess an equivalent
norm which is UKK ?}\\
Spaces with such a property are called UKK-able.  Known examples of
UKK-able spaces are the Hardy space $H_1$ of analytic functions on the
ball or on the polydisc \cite{bib:B.D.D.L}, the Lorentz spaces
$L_{p.1}(\mu)$
\cite{bib:C.D.L.T}, 
the trace class $C_1$ \cite{bib:L}, Gower's space which does not contain
$c_0, \,\, l_1,$ or any reflexive subspace \cite{bib:G}.
  Using  Prus ideas \cite{bib:Pr},
UKK-ability was intensively  studied in \cite{bib:D.G.K} for Banach spaces with a basis.\\

R. Huff associated an index to UKK-ability and proved that UKK-ability
implies that the index is finite.  He asked whether the converse is
true.  To our knowledge, it is still unknown in the general case.
G. Lancien \cite{bib:La} proved that the converse is true for the space of Bochner
integrable functions $L^p(X)$. More recently,  H. Knaust and E. Odell
proved that the converse is true for reflexive  Banach spaces
with  a finite dimensional decomposition  \cite{bib:K.O}.\\

In \cite{bib:K.T} a natural topology $\tau$ was introduced in Banach lattices.
This topology reduces to the topology of coordinatewise convergence in the case of
Banach spaces with basis, and to  the topology of convergence in measure in
usual spaces of measurable functions with order continuous norm.  In
this work, we define and investigate the UKK-ability for the topology
$\tau$. We first establish a generalization of a result of 
Kutzarova and Zachariades \cite{bib:K.Z} to the setting of
Banach lattices using well known
techniques of Pisier \cite{bib:Pi}. If $X$ is a Banach lattice
that satisfies a lower 1-estimate for any two disjoint
elements with a constant less than two, then $X$ must satisfy a lower
p-estimate for some $1\leq p <\infty.$ This will be Theorem 3.
 Our main result -Theorem 4- is that Banach lattices with a lower 
p-estimate for $1 \leq p < \infty$ are $\tau$-UKK-able.
 In particular, given that the $/tau$ topology and the topology of
convergence in measure are the same in $L_1$, we get the following 
result  
of \cite{bib:L2} that    $L_1$ is UKK  for the convergence in measure
  eventhough it fails
to be UKK-able for the weak-topology \cite{bib:D.G.K}.
 Let us also mention that $c_0$ is not $\tau$ UKK-able since it has an unconditional basis and it is known that it is not UKK-able
 for the weak topology, its Huff index being infinite. 

We also investigate the relationship between super-reflexivity
and super-UKK-ability and prove that these notions are equivalent in the 
lattice case. 

\section{Basic Definitions and Properties.} 

A sequence $(x_n)$ in a
Banach space $X$ is said to be $\varepsilon$-separated (with
$\varepsilon > 0$) if $$\inf \{ ||x_n - x_m|| ; \; n \not= m \} \geq
\varepsilon \ .$$ Let $\theta$ be any linear topology on $X$.  We say
that $X$ has the $\theta$-uniform Kadec-Klee (in short $\theta$-UKK)
property if for each $\varepsilon > 0$ there exists $\delta > 0$ such
that for  every $\varepsilon$-separated sequence $(x_n)$ in the unit ball
of $X$ which is $\theta$-convergent to $x$, we have that $||x|| \leq 1 -
\delta$.\\

 Let $X$ be a Banach lattice.  Let $u$ and $x$ be elements of $X$.
define $$ S_u(x) = x^+ \wedge |u| - x^- \wedge |u| $$ where $x^+$ and
$x^-$ are respectively the positive and negative parts of $x$ under
the lattice structure of $X$.\\

\bdfn
Let $X$ be a Banach lattice. The topology $\tau$ on $X$ is the coarsest linear
topology  for which the maps $x \rightarrow |||x|\wedge|u|||$
are continuous at $0$ for every $u \in X$.\\
\edfn

Thus a basis of neighborhoods of a point $x \in X$ is formed
 by the sets: $\{y \in X s.t. \| |y-x| \wedge u\| \leq \ep\}$
where $u \in X^+$ and $\ep > 0.$
The topology $\tau$ was introduced by M.A. Khamsi and Ph. Turpin in \cite{bib:K.T}
  to generalize a
fixed point theorem of P.K. Lin \cite{bib:Li} to the setting 
of Banach lattices.   It
was then observed that in the  spaces of measurable functions with order
continuous norm, $\tau$ is the topology of convergence in measure on
every set with finite measure.  On the other hand, for Banach spaces
with an unconditional basis, $\tau$ is the topology of coordinatewise
convergence. Among the nice properties of this topology is that it is
a Hausdorff, linear topology, coarser than  the norm; and if $K$
is $\tau$-compact then every sequence $(x_n)$ of points of $K$ has 
a $\tau$-convergent subsequence \cite{bib:K.T}.  \\

\begin{prop}
Let $X$ be a Banach lattice. For any $x$ and $y$ in $X$,
$x - S_y(x)$ and $y - S_x(y)$ are disjoints.
\end{prop}

\pf
  We have 
$$y - S_x(y) = y^+ - y^- - \Big[ y^+ \wedge |x| - y^- \wedge |x| \Big]
\ .$$
Therefore, we get
$$(1)\;\;\;\; \left\{ \begin{array}{lcl}
y - S_x(y) &=& (y^+ -  y^+ \wedge |x|) - (y^- - y^- \wedge |x|) \\
x - S_y(x) &=& (x^+ -  x^+ \wedge |y|) - (x^- - x^- \wedge |y|) 
\end{array} \right.$$

Note that $y^+ -  y^+ \wedge |x|\geq 0$ and $ y^- - y^- \wedge |x|
\geq 0$.  The same is true if we exchange $x$ with $y$.  Since $$y^+
\geq y^+ - y^+ \wedge |x| \;\;\;\mbox{and}\;\;\; y^- \geq y^- - y^-
\wedge |x| \ ,$$
we have
$$0 = y^+ \wedge y^- \geq \Big(y^+ - y^+ \wedge |x|\Big) \wedge \Big(
 y^- - y^- \wedge |x|\Big) \ .$$
The same is true if we exchange $y$ with $x$.  Thus the decomposition
in $(1)$ is a disjoint one in both equations.  In order to complete
the proof, we need to show that
$$(i)\;\;\; \Big(y^+ - y^+ \wedge |x|\Big) \wedge \Big(
 x^+ - x^+ \wedge |y|\Big) = 0$$
 and
$$(ii)\;\;\; \Big(y^+ - y^+ \wedge |x|\Big) \wedge \Big(
 x^- - x^- \wedge |y|\Big) = 0 \ .$$
Put $a = \Big(y^+ - y^+ \wedge |x|\Big) \wedge \Big(
 x^+ - x^+ \wedge |y|\Big)$.  It is clear that   
$0 \leq a \leq y^+ - y^+ \wedge |x|$; therefore, 
$$a + y^+ \wedge x^+ \leq a + y^+ \wedge |x| \leq y^+$$
and 
$$a + y^+ \wedge x^+ \leq a + |y| \wedge x^+ \leq x^+ \ .$$
Thus, $a + y^+ \wedge x^+ \leq y^+ \wedge x^+,$ which
implies that  $a \leq 0$.  This shows that we must have $a = 0$.  The proof of $(ii)$ is
done similarly.  
\qed

\begin{prop}
 Let $(x_n)$ be a
sequence $\tau$-convergent to $x$.  Then we have
$$\lim_{n \rightarrow \infty} ||S_u(x_n - x)|| = 0 \;\;\;\mbox{and}\;\;\;
\lim_{n \rightarrow \infty} ||S_{x_n - x}(u)|| = 0$$
for every $u$ in $X$.
\end{prop}

\pf
 By the triangle inequality we  have  $$||S_u(x_n - x)|| \leq
||(x_n - x)^+ \wedge |u| || + ||(x_n - x)^- \wedge |u| || \ .$$ Since
$X$ is a lattice, there exists $K >0$ such that $|x| \leq |y|$ implies
$||x|| \leq K ||y||$ for every $x$ and $y$ in $X$.  We then get
$$||S_u(x_n - x)|| \leq 2 K |||x_n - x| \wedge |u| ||$$ which implies that
$$\lim_{n \rightarrow \infty} ||S_u(x_n - x)|| = 0 \ .$$ On the other
hand, since $$||S_{x_n - x}(u)|| = |||x_n - x| \wedge u^+  -
|x_n - x| \wedge u^- || \leq |||x_n - x| \wedge u^+ || + ||
|x_n - x| \wedge u^- ||  $$
which clearly implies that
$$\lim_{n \rightarrow \infty} ||S_{x_n - x}(u)|| = 0 \ .$$
\qed

We now recall the definition of a lower p-estimate

\bdfn
 Let $1 \leq p < \infty$.  A
Banach lattice $X$ is said to satisfy a lower $p$-estimate (for
disjoint elements) if there exists a constant $C<\infty$ such that,
for every choice of pairwise disjoint elements $(x_i)_{1 \leq i \leq
n}$ in $X$, we have
$$\Bigg(\sum_{i} ||x_i||^p \Bigg)^{1/p} \leq C ||\sum_{i} x_i || \ .$$
\edfn 

\rmk {\bf 1:}
Let $X$ be a Banach space with an unconditional basis.
If $X$ is $\tau$-UKK-able then $X$ is UKK-able.

Indeed, as observed in \cite{bib:K.T}, in such a situation the $\tau$
topology
and the coordinatwise topology are the same. Given that a weakly
convergent sequence $(x_n)$ to an element $x \in X$ is coordinatwise
convergent to $x,$ the $\tau$ UKK-ability is clearly seen to be stronger 
than the UKK-ability.

\rmk {\bf 2:}
 Let $X$ be a Banach
lattice. 
Then $X$ admits an
equivalent norm $||.||_0$ which satisfies
$$(*)\;\;\;||x||_0 + ||y||_0 \leq 2 ||x+y||_0$$
for any disjoint $x,y \in X$.  

 Indeed, set
$$||x||_0 = \max\{||x^+||, ||x^-||\} \;\;\;\forall x \in X \ .$$
Then we have for any disjoint elements $x$, $y$ in $X$,
$$||x||_0 \leq ||x + y||_0 \;\;\;\mbox{and}\;\; ||y||_0 \leq ||x +
y||_0$$
which implies  $||x||_0 + ||y||_0 \leq 2 ||x + y||_0$.\\
We easily  check that $\| . \|_0 $ is an equivalent norm, since 
 $\|x\| = \|x^+ + x^-\| \leq \|x^+ \| + \|x^-\| \leq 2 \|x\|_0$
and $\|x\|_0 \leq ||x^+|| + ||x^-|| \leq  2 ||x||.$ \\

Note that every Banach lattice has an equivalent norm which satisfies a lower 1-estimate for any two disjoint elements with a constant $c \geq 2$.
 We don't have to consider the case $c=2$ since $c_0$
satisfies such an estimate with $c=2$ but is not UKK-able, and
therefore is not $\tau$-UKK-able by our first remark.
But for the case $c < 2,$ we do have a positive answer 
to the UKK-ability problem.
This will be our next theorem. 

We now  want to establish a generalization of a theorem by
Kutzarova and Zachariades \cite{bib:K.Z} and we will use techniques 
of Pisier \cite{bib:Pi}.
In the case of Banach spaces with basis this version was used by S.
Dilworth,
M. Girardi and D. Kutzarova in \cite{bib:D.G.K} as 
a criterion to show that Gower's space, that does not contain $c_0,
\,\, l_1$ or any reflexive subspace, is UKK-able. The proof is classical and we include it here for the sake of completeness.

\begin{thm}
Let $X$ be a Banach lattice. Suppose that there exists a constant $c,$
with $0<c<2,$ such that
$$\,\,(***)\,\,\,\, c \,
 \|x + y \| \geq (\|x\| + \|y\|)$$
for all disjoint $x$ and $y$ in $X.$
Then $X$ satisfies a lower p-estimate for some $1 \leq p <\infty.$
\end{thm}
  
\pf 
Suppose that there exists such a constant such that (***) is satisfied
for any disjoint $x_1$ and $x_2$ in $X.$
Then $\inf \Bigg(||x_1||, ||x_2||\Bigg) \leq (c/2) ||x_1 + x_2||.$\\
Similarly,
$\inf \Bigg( ||x_1||, ||x_2||, ||x_3||, ||x_4|| \Bigg) \leq (c/2)^2
||\sum_{i=1}^{4} x_i||$
for all $(x_i)_{1 \leq i \leq 4}$ that are disjoint.

So $$\inf_{1 \leq i \leq 2^n} \Bigg( ||x_i|| \Bigg) \leq (c/2)^n
||\sum_{i=1}^{2^n} x_i ||$$ for all $(x_i)_{1 \leq i \leq 2^n}$ that
are disjoint.

Now for a general $m,$ there exists a $k$ such that $2^k \leq m \leq
2^{k+1}.$
Then 
\begin{eqnarray*}
||\sum_{i=1}^{m} x_i || & = & || \sum_{i=1}^{2^k} x_i + \sum_{i=2^{k}
+1}^{m} x_i|| \\
                        & \geq & 1/c \Bigg(|| \sum_{i=1}^{2^k} x_i|| +|| \sum_{i=2^{k}
+1}^{m} x_i|| \Bigg) \\ 
                        & \geq & 1/c \,\,|| \sum_{i=1}^{2^k} x_i|| \\
                        & \geq & 1/c \,\,(2/c)^k \inf_{1 \leq i \leq 2^k}
			||x_i|| \\
                        & \geq & 1/c \,\,(2/c)^k \inf_{1 \leq i \leq m}
			||x_i|| \\   
\end{eqnarray*}

Thus $$\inf_{1 \leq i \leq m} ||x_i|| \leq \frac{c^{k+1}}{2^k} ||
\sum_{i=1}^{m} x_i||$$ for all finite disjoint sequences in $X.$

\noindent{\bf Claim 1:}{\em There exists a constant $c$ and a real
number $p< \infty$ such that 
$$\inf_{1 \leq i \leq m} ||x_i|| \leq \frac{c}{m^{1/p}}
||\sum_{i=1}^{m}x_i||$$
for all $m$ and disjoint $(x_i)_{1\leq i \leq m}$ in $X.$}
\\

\noindent{\bf Proof of Claim 1:}
It suffices then to prove that $(c/2)^k \leq
\frac{1}{m^{1/p}},$
which reduces to find a real number $p$ such that
$\frac{\log m}{k \log (2/c)} \leq p.$
But $2^k \leq m < 2^{k+1},$ so $2 \log 2 \geq (\frac{k+1}{k})\log 2
\geq \frac{\log m}{k}.$ Therefore $p = 2\frac{\log 2}{\log(2/c)}$
works, and the proof of Claim 1 is complete.
\\

\noindent{\bf Claim 2:}{\em For all reals $r > p,$ there exists a constant
$K_r$ so that for all $(x_i)$ that are disjoint in $X,$ we have
$$\Big(\sum_{i=1}^{n} ||x_i||^r \Big)^{1/r} \leq K_r 
|| \sum_{i=1}^{n}x_i||.$$}                  
\\

\noindent{\bf Proof of Claim 2:}
Suppose that $||x_1|| \leq ||x_2|| \leq ....\leq ||x_n||.$\\
Then $||x_1|| \leq \frac{c}{n^{1/p}} ||\sum_{i=1}^{n}x_i||,$ by Claim
1.

\noindent{Similarly,}
 $||x_2|| \leq \frac{c}{(n-1)^{1/p}} ||\sum_{i=2}^{n} x_i|| \leq c
\cdot  \frac{c}{(n-1)^{1/p}} ||\sum_{i=1}^{n} x_i||.$

\noindent{So} $||x_3|| \leq \frac{c}{(n-2)^{1/p}} ||\sum_{i=3}^{n} x_i|| \leq c
\cdot  \frac{c}{(n-2)^{1/p}} ||\sum_{i=1}^{n} x_i||,$

\noindent{and} $||x_n|| \leq \frac{c}{(n-(n-1))^{1/p}} ||x_{n}|| \leq c
\cdot  \frac{c}{(n-2)^{1/p}} ||\sum_{i=1}^{n} x_i||.$

\noindent{Combining all inequalities together we get:}
$$\sum_{i=1}^{n}||x_i||^r \leq
\Big(\sum_{i=1}^{n}\frac{c^{2r}}{i^{r/p}}\Big)||\sum_{i=1}^{n}x_i||^r.$$
So it suffices to let $K_r = \Big(\sum_{i=1}^{\infty}
\frac{c^{2r}}{i^{r/p}}\Big)^{1/r}.$
The proof of Claim 2 is complete, and the theorem is proved.
\qed

\section{Main result.}

As in \cite{bib:H}, assume that $X$ is a Banach lattice satisfying a lower
$p$-estimate ($1 \leq p < \infty$), and  set 
\begin{equation}
||x||_0 = \sup\Bigg(\sum_{n} ||x_n||^p \Bigg)^{1/p} \;,
\label{norm0}
\end{equation}
 where the supremum is
taken over all pairwise disjoint elements $(x_i)$ such that $x =
\sum_{i} x_i$.  It is easy to see that $||.||_0$ is a
norm on $X$ for which $||.|| \leq ||.||_0 \leq C ||.||$, i.e. the
norm $||.||_0$ is equivalent to $||.||$.  Moreover we have
$$\;\;\;\; ||x||_0^p + ||y||_0^p \leq ||x + y||_0^p $$
for any disjoint elements $x$ and $y$ in $X$.\\

\begin{thm}
   Let $X$ be a Banach
lattice with a $p$-lower estimate for $ 1 \leq p < \infty$.  Then the
norm $||.||_0$ (defined by (1)) is $\tau$-UKK.
\end{thm}

\pf
  Let $\varepsilon > 0$ and $(x_n)$ be
$\varepsilon$-separated elements in the unit ball of $X$  such that
$\tau-\lim x_n = x$.  Since for every integers $n \not= m$ we have
$$\varepsilon \leq ||x_n - x_m ||_0 \leq ||x_n - x||_0 + ||x_m - x||_0$$
we get
$$\frac{\varepsilon}{2} \leq \liminf_{n \rightarrow \infty} ||x_n -
x||_0 \ .$$
Since $x_n - x - S_{|x|}(x_n - x)$ and $x - S_{|x_n-x|}(x)$ are
disjoints (by Proposition 1), then we have
$$||x_n - x - S_{|x|}(x_n - x)||_0^p + ||x - S_{|x_n-x|}(x)||_0^p \leq
||x_n - x - S_{|x|}(x_n - x) + x - S_{|x_n-x|}(x)||_0^p \ .$$
But $x_n - x - S_{|x|}(x_n - x) + x - S_{|x_n-x|}(x) = 
x_n  - S_{|x|}(x_n - x)  - S_{|x_n-x|}(x)$, and since
$$\lim_{n \rightarrow \infty}||S_{|x|}(x_n - x)||_0 = 0 = 
\lim_{n \rightarrow \infty} || S_{|x_n-x|}(x)||_0,$$
then we have
\begin{eqnarray*}
\liminf_{n \rightarrow \infty}||x_n-x||^p_0 + ||x||_0^p & = &
\liminf_{n \rightarrow \infty} ||x_n - x - S_{|x|}(x_n - x)||_0^p + 
\liminf_{n \rightarrow \infty}||x - S_{|x_n-x|}(x)||_0^p \\
                 &  \leq &
\liminf_{n \rightarrow \infty}||x_n||^p_0 \leq 1. 
\end{eqnarray*}
Therefore, we obtain
$$||x||^p_0 \leq 1 - \liminf_{n \rightarrow \infty}||x_n-x||^p_0 \leq
1 - \left(\frac{\varepsilon}{2}\right)^p $$
which implies that $||x||_0 \leq 1-\delta$, with
$$\delta = 1 -\left(1 -
\left(\frac{\varepsilon}{2}\right)^p\right)^{1/p} \ .$$
This completes the proof of Theorem 2.\\

We then get the following 

\begin{cor}
 Let $X$ be a Banach
lattice that does not contain $l_{\infty}^{n}$ uniformly for every
$n.$
  Then $X$ 
is $\tau$-UKK-able.
\end{cor}

\pf
Let $X$ be a Banach lattice which does not contain $\l_{\infty}^{n}$
uniformly for all $n.$ By the result of \cite{bib:J} and \cite{bib:Sh}
$X$ has a lower p-estimate for some $1 \leq p < \infty.$
We then use the previous result.
\qed

\section{Super-UKK-ability in Banach lattices}
In \cite{bib:D.G.K} the following question is raised:
Is
super-reflexivity and super-UKK-ability equivalent in any Banach
space? In the next proposition, we give a positive answer in  the case
of Banach lattices. To our knowledge the general question is still open.

\begin{prop}
Let $X$ be a Banach lattice. $X$ is super-reflexive if and only if
$X$ is super-UKK-able.
\end{prop}

\pf
If $X$ is super-reflexive then $X$ is clearly super-UKK-able.
Suppose now that $X$ is super-UKK-able. 
Take any ultrapower $\Pi X/{\cal U}.$ 
This is still a Banach lattice.
So it either contain $c_0,$
or $l_1,$ or is reflexive. 
It can't contain $c_0$ since $\Pi X/{\cal U}$ is UKK-able.
If $\Pi X/{\cal U}$ contains $l_1$ then its ultrapower
$\Pi\Bigg(\Pi X/{\cal U}\Bigg)$ would contain $L_1$ 
which is not UKK-able.
Therefore $\Pi X/{\cal U}$ is reflexive, i.e. $X$ is super-reflexive.
\qed    

\rmk

As it was pointed out to us by D. Leung, the proof of the previous
proposition
actually shows that if $X$ is a super-UKK-able Banach space,
then $X$ cannot contain $l_{\infty}^n$ or $l_1^n$ uniformly in $n.$
However, it is well known that in 
a general Banach space, the latter condition is not sufficient
for super-reflexivity. In fact G. Pisier and Q. Xu constructed for
every $q>2$ a non-reflexive space of type 2 and cotype $q.$

\bigskip

Department of Mathematical Sciences

The University of Missouri at Columbia

65211 Columbia, Missouri.

\&

Department of Mathematical Sciences 

The University of Texas at El Paso

El Paso, Texas 79968-0514.

\end{document}